\documentclass[12pt]{amsart}
\usepackage{latexsym, amssymb, amsmath, amscd}
 \usepackage{eucal}

    \newtheorem{rema}{Remark}[section]
    \newtheorem{propo}[rema]{Proposition}

   \newtheorem{theo}[rema]{Theorem}
   \newtheorem{def-theo}[rema]{Definition-Theorem}
 \newtheorem{conj}[rema]{Conjecture}
   \newtheorem{defi}[rema]{Definition}
    \newtheorem{lemma}[rema]{Lemma}
    \newtheorem{corol}[rema]{Corollary}
     
  \newtheorem{rmk}[rema]{Remark}

	\newcommand{\nno}{\nonumber}

 \newcommand{\pf}{{\it Proof:}\hspace{2ex}}
 
 \newcommand{\epfv}{\hspace{1em}$\Box$\vspace{1em}}

\newcommand{\bZ}{{\mathbb Z}}

\newcommand{\bN}{{\mathbb N}}

\newcommand{\cA}{{\mathcal A}}

\newcommand{\I}{{\operatorname I}}

\newcommand{\Ad}{\mbox{\rm Ad\,}}
\newcommand{\ad}{\mbox{\rm ad\,}}

\newcommand{\cB}{{\mathcal B}}

\newcommand{\cE}{{\mathcal E}}
\newcommand{\cD}{{\mathcal D}}

\newcommand{\Ker}{\mtype{\rm Ker\,}}
\newcommand{\im}{\mtype{Im\,}}

\newcommand{\mtype}{\operatorname}

\newcommand{\Der}{{\cD er}}
\newcommand{\Eder}{{\cE der}}
\newcommand{\cEnd}{{\cE nd}}

\title[The LNED and LFED Conjectures for Algebraic Algebras]
{The LNED and LFED Conjectures for Algebraic Algebras}

  \author{Wenhua Zhao}      
    \date{\today}
\address{Department of Mathematics, Illinois State University, Normal, IL 61761. Email: wzhao@ilstu.edu}

\begin{document}

\begin{abstract}
Let $K$ be a field of characteristic zero and $\cA$ a $K$-algebra such that all 
the $K$-subalgebras generated by finitely many elements of $\cA$ are finite dimensional over $K$. A $K$-$\cE$-derivation of $\cA$ is a $K$-linear map of the form $\I-\phi$ for 
some $K$-algebra endomorphism $\phi$ of $\cA$, where $\I$ denotes the identity map of $\cA$. 
In this paper we first show that for all locally finite $K$-derivations $D$ and locally finite $K$-algebra automorphisms $\phi$ of $\cA$, the images of $D$ and $\I-\phi$ do not contain any nonzero idempotent of $\cA$. We then use this result to show some 
cases of the LFED and LNED conjectures proposed in \cite{Open-LFNED}. 
More precisely, 
We show the LNED conjecture for $\cA$, and the LFED conjecture for 
all locally finite $K$-derivations of $\cA$ and all locally finite $K$-$\cE$-derivations of  the form $\delta=\I-\phi$ with $\phi$ being surjective. In particular,  
both conjectures are proved for all finite dimensional $K$-algebras. 
Furthermore, some finite extensions of derivations and automorphism 
to inner derivations and inner automorphisms, respectively, have also been established. 
This result is not only crucial in the proofs of the results above, but also interesting 
on its own right. 
\end{abstract}

\keywords{Mathieu subspaces (Mathieu-Zhao spaces), the LNED conjecture, the LFED conjecture,   locally finite or locally nilpotent derivations and $\cE$-derivations, inner derivations, inner automorphisms, idempotents}
   
\subjclass[2000]{47B47, 08A35, 16W25, 16D99}








%
%



\thanks{The author has been partially supported 
by the Simons Foundation grant 278638}

 \bibliographystyle{alpha}
    \maketitle


\renewcommand{\theequation}{\thesection.\arabic{equation}}
\renewcommand{\therema}{\thesection.\arabic{rema}}
\setcounter{equation}{0}
\setcounter{rema}{0}
\setcounter{section}{0}

\section{\bf Motivations and the Main Results} \label{S1}

Let $R$ be a unital ring (not necessarily commutative) and $\cA$ an $R$-algebra. 
We denote by $1_\cA$ or simply $1$ the identity element of $\cA$, if $\cA$ is unital, and 
$\I_\cA$ or simply $\I$ the identity map of $\cA$, if $\cA$ is clear in the context.  

An $R$-linear endomorphism $\eta$ of $\cA$ is said to be {\it locally nilpotent} (LN) 
if for each $a\in \cA$ there exists $m\ge 1$ such that $\eta^m(a)=0$, 
and {\it locally finite} (LF) if for each $a\in \cA$ the $R$-submodule spanned 
by $\eta^i(a)$ $(i\ge 0)$ over $R$ is finitely generated.    

By an $R$-derivation $D$ of $\cA$ we mean an $R$-linear map 
$D:\cA \to \cA$ that satisfies $D(ab)=D(a)b+aD(b)$ for all $a, b\in \cA$. 
By an $R$-$\cE$-derivation $\delta$ of $\cA$ 
we mean an $R$-linear map $\delta:\cA \to \cA$ such that for all 
$a, b\in \cA$ the following equation holds:
\begin{align}\label{ProdRule2}
\delta(ab)=\delta(a)b+a\delta(b)-\delta(a)\delta(b).  
\end{align}

It is easy to verify that $\delta$ is an $R$-$\cE$-derivation of $\cA$, 
if and only if $\delta=\I-\phi$ for some $R$-algebra endomorphism 
$\phi$ of $\cA$. Therefore an $R$-$\cE$-derivation is  
a special so-called $(s_1, s_2)$-derivation introduced 
by N. Jacobson \cite{J} and also a special 
semi-derivation introduced 
by J. Bergen in \cite{Bergen}. $R$-$\cE$-derivations have also been 
studied by many others under some different names such as 
$f$-derivations in \cite{E0, E} and 
$\phi$-derivations in \cite{BFF, BV}, etc..

We denote by $\cEnd_R(\cA)$ the set of all 
$R$-algebra endomorphisms of $\cA$, $\Der_R(\cA)$ the set of all 
$R$-derivations of $\cA$, and $\Eder_R(\cA)$ the set of all 
$R$-$\cE$-derivations of $\cA$. Furthermore, for each $R$-linear endomorphism 
$\eta$ of $\cA$ we denote by $\im \eta$ the {\it image} of $\eta$, i.e., 
$\im\eta\!:=\eta(\cA)$, and $\Ker \eta$ the {\it kernel} of $\eta$.
When $\eta$ is an $R$-derivation or $R$-$\cE$-derivation, we also denote  
by $\cA^{{}^\eta}$ the kernel of $\eta$.  
 
Next let us recall the following notion first introduced in \cite{GIC, MS}.

\begin{defi} \label{Def-MS}
Let $\vartheta$ represent the words:  
$\mtype{left}$, $\mtype{right}$, or $\mtype{two-sided}$. An $R$-subspace $V$ 
of an $R$-algebra $\cA$ is said to be a $\vartheta$-Mathieu subspace ($\vartheta$-MS) 
of $\cA$ if for all $a, b, c\in \cA$ with $a^m\in V$ for all $m\ge 1$,  
the following conditions hold: 
\begin{enumerate}
  \item[$1)$] $ba^m\in V$ for all $m\gg 0$, if $\vartheta=\mtype{left}$;
  \item[$2)$] $a^mc\in V$ for all $m\gg 0$, if $\vartheta=\mtype{right}$;
  \item[$3)$] $b a^m c \in V$ for all $m\gg 0$, if $\vartheta=\mtype{two-sided}$. 
\end{enumerate} 
\end{defi}

A two-sided MS will also be simply called a MS. 
Note that a MS is also called a {\it Mathieu-Zhao space} 
in the literature (e.g., see \cite{DEZ, EN, EH}, etc.) 
as first suggested by A. van den Essen \cite{E2}. 
 
The introduction of the new notion 
is mainly motivated by the study in \cite{Ma, IC} of the well-known 
Jacobian conjecture (see \cite{K, BCW, E}). 
See also \cite{DEZ}. But, a more interesting aspect 
of the notion is that it provides a natural but highly non-trivial   
generalization of the notion of ideals of associative algebras. 

Note that the MSs of algebraic algebras over a field $K$ can be characterized 
by the following theorem, which is a special case of \cite[Theorem $4.2$]{MS}.

\begin{theo} \label{IdemCri} 
Let $K$ be a field of arbitrary characteristic 
and $\cA$ a unital $K$-algebra that is algebraic over $K$. 
Then a $K$-subspace $V$ is a MS of $\cA$, 
if and only if for every idempotent $e$ (i.e., $e^2=e$), the principal 
ideal $(e)$ of $\cA$ generated by $e$ is contained in $V$.
\end{theo} 

Next let us recall the following two conjectures proposed in \cite{Open-LFNED}. 

\begin{conj}[{\bf The LFED Conjecture}]\label{LFED-Conj}
Let $K$ be a field of characteristic zero, $\cA$ a $K$-algebra 
and $\delta$ a LF (locally finite) $K$-derivation or a LF 
$K$-$\cE$-derivation of $\cA$. Then the image 
$\im\delta$ of $\delta$ is a MS of $\cA$.   
\end{conj}

\begin{conj}[{\bf The LNED Conjecture}]\label{LNED-Conj}
Let $K$ be a field of characteristic zero, $\cA$ a $K$-algebra 
and $\delta$ a LN (locally nilpotent) $K$-derivation or a LN 
$K$-$\cE$-derivation of $\cA$. Then $\delta$ maps every 
$\vartheta$-ideal of $\cA$ to a $\vartheta$-MS of $\cA$, 
where $\vartheta$ represents the words:  
$\mtype{left}$, $\mtype{right}$, or $\mtype{two-sided}$. 
\end{conj}
 
In this paper we prove some cases of the LFED and LNED 
conjectures above under the following condition on $\cA$:  

\begin{enumerate}
\item[$(KC)$] \label{KC}  
{\it all the $K$-subalgebras generated by finitely 
many elements of $\cA$ are finite dimensional over $K$.}   
\end{enumerate}

For the studies of some other cases of the LFED and LNED conjectures above, 
see \cite{EWZ}, \cite{Open-LFNED}--\cite{OneVariableCase}.

\begin{rmk}
The condition $(KC)$ above is satisfied by all commutative algebraic 
$K$-algebras. But for noncommutative algebraic $K$-algebras the condition $(KC)$ is the same as saying that the well-known Kurosch's Problem \cite{Kurosch} (see also \cite{R} and \cite{Zel}) has a positive answer for all finitely generated $K$-subalgebras of $\cA$. Note that the Kurosch's Problem does not have a positive answer for all noncommutative affine algebras. 
\end{rmk}

In this paper we shall first show  
the following   
 
\begin{theo} \label{Idem-Alg-LNED}
Let $K$ be a field of characteristic zero and $\cA$ a 
$K$-algebra that satisfies the condition $(KC)$ above. 
Let $\delta$ be a LF $K$-derivation of 
$\cA$ or a LF $K$-$\cE$-derivation of $\cA$. Then the image  
$\im\delta$ of $\delta$ does not contain 
any non-zero idempotent $e$ of $\cA$ if $\delta$ also satisfies 
one of the following two conditions: 
\begin{enumerate}
\item[$1)$]  $\delta\in \Der_K(\cA)$; 
\item[$2)$]  $\delta=\I-\phi$ for some    
$K$-algebra automorphism $\phi$ of $\cA$.
\end{enumerate}
 \end{theo}

Note that for every LN $K$-$\cE$-derivation
$\delta=\I-\phi$, the $K$-algebra endomorphism $\phi$ 
is invertible with the inverse map 
$(\I-\delta)^{-1}=\sum_{i=0}^\infty \delta^i$. 
Note also that every LN $K$-linear map is also LF. Therefore,     
by Theorems \ref{IdemCri} and \ref{Idem-Alg-LNED} above 
we have the following  

\begin{corol}\label{Alg-LNED}
Let $K$, $\cA$ be as in Theorem \ref{Idem-Alg-LNED} and $\delta\in \Der_K(\cA)$ or 
$\Eder_K(\cA)$. Then the following statements hold:
\begin{enumerate}
  \item[$1)$] if $\delta$ is LN, then $\delta$ maps every $K$-subspace 
of $\cA$ to a MS of $\cA$. In particular, the LNED conjecture \ref{LNED-Conj} holds for $\cA$; 
  \item[$2)$] if $\delta$ is LF and satisfies the condition $1)$ or $2)$ in   
Theorem \ref{Idem-Alg-LNED}, then $\im \delta$ is a MS of $\cA$, i.e.,  
the LFED conjecture \ref{LFED-Conj} holds for $\delta$.
\end{enumerate}
\end{corol}

%
%

Actually, it will be shown in Proposition \ref{LFED-Prop} in subsection \ref{S2.3} 
that the LFED conjecture \ref{LFED-Conj} holds also for all the LF 
$K$-$\cE$-derivations $\delta$ of $\cA$ of 
the form $\delta=\I-\phi$ with $\phi\in \cEnd_K(\cA)$ being surjective.
Some other cases of the LFED conjecture \ref{LFED-Conj} and 
the LNED conjecture \ref{LNED-Conj} will also be proved in 
subsection \ref{S2.3}. For example, it will be shown in 
Proposition \ref{FinDimCase} that both the LFED and LNED conjectures 
hold for all finite dimensional algebras 
over a field of characteristic zero.  

\begin{rmk}
It is easy to see that the following algebras  
over a field $K$ satisfy the condition $(KC)$:
\begin{enumerate}
  \item[$i)$] finite dimensional algebras over $K$;
  \item[$ii)$] commutative algebras that are algebraic over $K$;
  \item[$iii)$] the union of an increasing sequence of  
  $K$-algebras in $i)$ or $ii)$, e.g., the $\bN\times\bN$ matrix algebra 
over $K$ with only finitely nonzero entries; etc.    
\end{enumerate}  
Therefore, Theorem \ref{Idem-Alg-LNED} and Corollary \ref{Alg-LNED} 
above as well as some other results proved in this paper 
apply to all the algebras above.  
\end{rmk}

In order to show Theorem \ref{Idem-Alg-LNED} we need to show a  
theorem (Theorem \ref{FinExt} below) on finite extensions of a derivation  
(resp.,, an automorphism) of $\cA$ to an inner derivation 
(resp.,, an inner automorphism). Since this theorem is interesting on 
its own right, we formulate and prove  
it in a more general setting. 

Let $R$ be a unital ring (not necessarily commutative), 
$\cA$ an $R$-algebra. Recall that a derivation $D$ of 
$\cA$ is {\it inner} if there exists $u\in \cA$ such that $D=\ad_u$, 
where $\ad_u$ is the adjoint derivation induced by $u$, 
i.e., $\ad_u(a)\!:=[u, a]\!:=ua-au$ for all $a\in \cA$. 
An automorphism $\phi$ of $\cA$ is {\it inner} if 
there exists a unit $u\in \cA$ such that $\phi=\Ad_u$, 
where $\Ad_u$ is the conjugation automorphism of $\cA$ 
induced by $u$, i.e., $\Ad_u(a)\!:=uau^{-1}$ for 
all $a\in \cA$. 

An $R$-linear endomorphism $\eta$ 
of $\cA$ is {\it integral} over $R$ 
if there exists a monic polynomial 
$p(t)\in R[t]$ such that $p(\eta)=0$. 
Although $R$ may not be commutative, the valuation 
$p(\eta)$ does not depend on which side 
we write the coefficients of $p(t)$, since $\eta$ 
is an $R$-linear endomorphism of $\cA$. 
  
\begin{theo}\label{FinExt}
Let $R$ be a unital ring (not necessarily commutative), 
$\cA$ a unital $R$-algebra and $D$ (resp., $\phi$) an $R$-derivation (resp., $R$-algebra automorphism) of $\cA$. Let $p(t)$ be a monic polynomial 
in $R[t]$ such that $p(D)=0$ (resp., $p(\phi)=0$ and $p(0)$ is a unit of $R$). 
Then there exists an $R$-algebra 
$\cB$ containing $\cA$ as an $R$-subalgebra such that 
the following statements hold:
\begin{enumerate}
\item[$1)$] $\cB$ is finitely generated as 
both a left $\cA$-module and a right $\cA$-module;   
\item[$2)$] there exists $u\in \cB$ (resp., a unit $u\in \cB$) 
such that $p(u)=0$ and $D=\ad_u\, |_{{}_\cA}$ 
(resp., $\phi=\Ad_u  |_{{}_\cA}$).     
\end{enumerate}    
\end{theo}

{\bf Arrangement:}  We give a proof for the derivation case of Theorem \ref{Idem-Alg-LNED} in subsection \ref{S2.1}, and the $\cE$-derivation case 
in a more general setting in subsection \ref{S2.2} 
(see Proposition \ref{IdemImI-Phi}). 
In subsection \ref{S2.3} we discuss some other consequences of 
Theorem \ref{Idem-Alg-LNED} and Proposition \ref{IdemImI-Phi}.  
We then give a proof for the derivation case of Theorem \ref{FinExt} 
in subsection \ref{S3.1}, and the automorphism case in 
subsection \ref{S3.2}.  \\

{\bf Acknowledgment:} The author is very grateful to Professors Arno van de Essen for reading carefully an earlier version of the paper and pointing out some mistakes and typos, etc..

\renewcommand{\theequation}{\thesection.\arabic{equation}}
\renewcommand{\therema}{\thesection.\arabic{rema}}
\setcounter{equation}{0}
\setcounter{rema}{0}

\section{\bf Proof and Some Consequences of Theorem \ref{Idem-Alg-LNED}}\label{S2}

In this section we assume Theorem \ref{FinExt} and   
give a proof for Theorem \ref{Idem-Alg-LNED}. We divide the proof into two cases:  
the case of derivations in subsection \ref{S2.1} and the case 
of $\cE$-derivations in subsection \ref{S2.2}. We then derive some 
consequences of Theorem \ref{Idem-Alg-LNED} and Proposition \ref{IdemImI-Phi} 
in subsection \ref{S2.3}.

Throughout this section $K$ denotes a field of 
characteristic zero and {\bf $\cA$ a unital $K$-algebra 
that satisfies  the condition} $\mathbf{(KC)}$ on page \pageref{KC}.
All the notations introduced in Section \ref{S1} will   
also be freely used.

\subsection{Proof of Theorem \ref{Idem-Alg-LNED}, 1)}\label{S2.1}
Let $D$ be a LF (locally finite) $K$-derivation of $\cA$. 
Assume that Theorem \ref{Idem-Alg-LNED}, $1)$ fails, 
i.e., there exists a nonzero idempotent $e\in \im D$.  
Let $a\in \cA$ such that $e=D a$ and $\cA_1$ be the $K$-subalgebra 
generated by $D^i(a^j)$ $(i, j\ge 0)$ over $K$. 
Then $\cA_1$ is $D$-invariant. Furthermore, 
by the condition $(KC)$ assumed on $\cA$ $a$ is 
algebraic over $K$. Therefore there exists 
$d\ge 1$ (e.g., the degree of the minimal polynomial of $a$ over $K$) 
such that $a^j$ ($j\ge d$) lies in the 
$K$-subspace spanned over $K$ 
by $a^i$ $(0\le i\le d-1)$.

On the other hand, since $D$ is LF, for each fixed 
$0\le i\le d-1$ the $K$-subspace spanned over $K$ 
by $D^j(a^i)$ $(j\ge 0)$ is of finite dimension over $K$. 
Therefore $\cA_1$ is generated by finitely 
many elements of $\cA$. By the condition $(KC)$ 
on $\cA$ again $\cA_1$ is also finite 
dimensional over $K$. Replacing $\cA$ by 
$\cA_1$ and $D$ by $D\mid_{\cA_1}$ we may assume 
that $\cA$ itself is finite dimensional 
over $K$, and consequently, $D$ is a $K$-derivation 
of $\cA$ that is integral over $K$.  

By Theorem \ref{FinExt} there exists a $K$-algebra 
extension $\cB$ of $\cA$ and some $u\in \cB$ such that 
$D=\ad_u\,|_\cA$. Furthermore, since $\cB$ is finitely generated 
as a left $\cA$-module, $\cB$ is also 
finite dimensional over $K$.

Let $n=\dim_K \cB$ and $\mu$ be the regular representation 
of $\cB$ to the $K$-algebra  $\operatorname{Hom}_K(\cB, \cB)$  
of all $K$-linear endomorphisms of $\cA$, i.e., 
$\mu(b)$ for each $b \in \cB$ is the multiplication map 
by $b$ from the left. Then $\mu$ is a 
faithful representation, since $\cA$, hence also $\cB$, 
is unital. Choosing a $K$-linear basis of $\cB$ we identify 
$\operatorname{Hom}_K(\cB, \cB)$ with the $n\times n$ 
matrix algebra $M_n(K)$ over $K$. Since 
\begin{align}\label{S2.1-peq1}
\mu(e)=\mu(\ad_u(a))=\mu([u, a])=[\mu(u), \mu(a)], 
\end{align}  
we see that the trace of the matrix $\mu(e)$ is equal 
to zero. 

On the other hand, $\mu(e)$ is  
also a nonzero idempotent matrix. It is well-known 
in linear algebra that the trace of any nonzero idempotent 
in $M_n(K)$ is not zero. For example, by using the Jordan form 
of the matrix $\mu(e)$ in $M_n(\bar K)$, where $\bar K$ is the algebraic 
closure of $K$, it is easy to see that the trace of $\mu(e)$ is actually equal 
to its rank. Hence we get a contradiction. Therefore,  
statement $1)$ in Theorem \ref{Idem-Alg-LNED} holds.

\subsection{The $\cE$-Derivation Case of Theorem \ref{Idem-Alg-LNED}} \label{S2.2}

Throughout this and also the next subsection we fix a $K$-algebra 
endomorphism $\phi$ of $\cA$ and set 
$\Ker_{\ge1}\phi\!:=\sum_{i=1}^\infty \Ker \phi^i$ 
and $\bar \cA\!:=\cA/\Ker_{\ge1}\phi$. We denote by $\pi$ the quotient map from $\cA$ to $\bar\cA$ and $\bar\phi$ the induced map by $\phi$ from $\bar\cA$ to $\bar\cA$.

We will also freely use the fact that $\phi$ is LF (locally finite), 
if and only if $\I-\phi$ is LF.  
 
\begin{lemma} \label{Lma2.2-1} 
With the setting above we have  
\begin{enumerate}
  \item[$1)$] $\Ker_{\ge 1} \phi \subseteq \im(\I-\phi)$.
  \item[$2)$] $\bar\phi$ is injective. 
  \end{enumerate}
\end{lemma}

This lemma is actually part of \cite[Lemma 5.3]{Open-LFNED}. 
But for the sake of completeness we include a proof here. \\

\underline{\it Proof of Lemma \ref{Lma2.2-1}}: $1)$ Let $a\in \Ker_{\ge 1} \phi$. Then $\phi^k(a) = 0$ for some $k \ge 1$. 
Let $b=\sum_{i=0}^\infty \phi^i(a)$, which is a well-defined element of $\cA$. 
Then $(\I-\phi)(b)=(\I-\phi)(\sum_{i=1}^\infty \phi^i)(a)=a$. Therefore $a\in \im (\I-\phi)$. 
 
$2)$ Let $a\in \cA$ such that $\bar \phi(\pi(a))=0$. Since $\bar\phi\pi=\pi\phi$, 
we have $\pi(\phi(a))=0$, i.e., $\phi(a)\in \Ker_{\ge 1}\phi$. Then  
$\phi^{k+1}(a)=\phi^k(\phi(a))=0$ for some $k\ge 1$, and   
$a\in \Ker_{\ge 1}\phi=\Ker\pi$, whence $\pi(a)=0$ and $\bar\phi$ is injective.
\epfv



Next, we show the following proposition, from which Theorem \ref{Idem-Alg-LNED}, 
$2)$ follows immediately, since $\phi$ is LF, if  and only if $\I-\phi$ is LF.

\begin{propo}\label{IdemImI-Phi}
Assume that $\cA$ satisfies the condition $(KC)$, and 
$\phi\in\cEnd_K(\cA)$ such that $\phi$ is LF and 
$\bar\phi$ is surjective (e.g., when $\phi$ itself is surjective). 
Then for each idempotent $e\in \cA$ 
we have that $e\in \im (\I-\phi)$, if and only if  
$e\in \Ker_{\ge 1}\phi$. 

Consequently, if $\phi$ is a $K$-algebra automorphism of $\cA$, 
then  $\im (\I-\phi)$ does not contain 
any nonzero idempotent of $\cA$. 
\end{propo}
 
\pf We first consider the case that $\phi$ is bijective. 
Note that in this case $\Ker_{\ge 1}\phi=0$ and $\bar\cA=\cA$. 
So it suffices to show that $\im(\I-\phi)$ does not contain 
any nonzero idempotent of $\cA$. But this can be proved 
by a similar argument as the proof of 
Theorem \ref{Idem-Alg-LNED}, $1)$ 
in the previous subsection.  
For example, by letting $p(t)$ be the minimal 
polynomial of $\phi$ we have that $p(0)\ne 0$ and hence is 
a unit of the base ring $K$. So we may apply 
Theorem \ref{FinExt} to the automorphism $\phi$ 
(instead of applying it to the derivation $D$). 
By using the same notation Eq.\,(\ref{S2.1-peq1}) becomes 
\begin{align*}
\mu(e)=\mu(a-\Ad_u(a))=\mu(a-ua u^{-1})=\mu(a)-\mu(u)\mu(a) \mu(u)^{-1}, 
\end{align*} 
from which we see that the trace of the nonzero idempotent matrix 
$\mu(e)$ is equal to zero, which is a contradiction again. 
  
Next we consider the general case. Note that by  
Lemma \ref{Lma2.2-1}, $1)$ it suffices to show that 
every idempotent $e\in\im(\I_\cA-\phi)$ lies in 
$\Ker_{\ge 1}\phi$. 

By Lemma \ref{Lma2.2-1}, $2)$ and the surjectivity of $\bar \phi$ 
we see that $\bar \phi$ is a $K$-algebra automorphism of $\bar \cA$.  
Then by the bijective case shown above $\im(\I_{\bar\cA}-\bar\phi)$ 
does not contain any nonzero idempotent of $\bar\cA$.

On the other hand, for all $e\in\im(\I_\cA-\phi)$ it is easy to see that 
$\pi(e)$ is an idempotent in $\im(\I_{\bar\cA}-\bar\phi)$. Therefore we have 
$\pi(e)=0$, and hence $e\in \Ker\pi=\Ker_{\ge 1}\phi$, as desired. 
\epfv

%
%
%

\subsection{Some Consequences}\label{S2.3}  
In this subsection we derive some 
consequences of Theorem \ref{Idem-Alg-LNED} 
and Proposition \ref{IdemImI-Phi}. 
All the notations fixed in the previous subsection 
will still be in force in this subsection.

\begin{propo} \label{LFED-Prop}
Let $\delta$ be a LF $K$-derivation of $\cA$, or 
a LF $K$-$\cE$-derivation of the form 
$\delta=\I-\phi$ for some $\phi\in\cEnd_K(\cA)$ such that 
$\bar\phi$ is surjective (e.g., when $\phi$ itself is surjective). 
Then the LFED conjecture \ref{LFED-Conj}
holds for $\delta$.  
\end{propo} 

\pf If $\delta$ is a $K$-derivation of $\cA$,  
then the proposition follows immediately from 
Corollary \ref{Alg-LNED}, $2)$. 

Assume that $\delta=\I-\phi$ for some $\phi\in\cEnd_K(\cA)$ such that 
$\bar\phi$ is surjective. Let $e$ be 
an idempotent lying in $\im\delta$. Then $e\in \Ker_{\ge 1}\phi$ by 
Proposition \ref{IdemImI-Phi}. 
Since $\Ker_{\ge 1}\phi$ is an ideal of $\cA$, the principal ideal 
$(e)\subseteq \Ker_{\ge 1}\phi$. Then   
by Lemma \ref{Lma2.2-1}, $1)$ we have 
$(e)\subseteq \im(\I-\phi)$, and 
by Theorem \ref{IdemCri} the proposition 
follows.    
\epfv
  
%
%

\begin{propo}\label{FinDimCase}
Both the LFED conjecture \ref{LFED-Conj} and the LNED conjecture \ref{LNED-Conj} hold  
for all finite dimensional algebras over a field $K$ of characteristic zero.
\end{propo}

\pf Let $\cA$ be a finite dimensional $K$-algebra. Then by Corollary \ref{Alg-LNED} 
it is sufficient to show the LFED conjecture \ref{LFED-Conj} for every  
LF $K$-$\cE$-derivation $\delta$ of $\cA$. 

Write $\delta=\I-\phi$ for some $\phi\in\cEnd_K(\cA)$.  
Then by Lemma \ref{Lma2.2-1}, $1)$ we have 
that $\bar\phi$ is injective. Since $\cA$ is finite dimensional 
over $K$, then so is $\bar\cA$. Hence $\bar\phi$ is also surjective. Then   
by Proposition \ref{LFED-Prop} the LFED conjecture \ref{LFED-Conj} 
holds for $\delta$, whence the proposition follows.   
\epfv

Finally, it is also worthy to point out the following special case of 
the LFED conjecture \ref{LFED-Conj} and the LNED conjecture \ref{LNED-Conj}, 
which follows directly from Theorems \ref{IdemCri}, \ref{Idem-Alg-LNED} and 
the fact that all finite order $K$-algebra automorphisms of $\cA$ are LF.

\begin{corol}
For every finite order $K$-algebra automorphism $\phi$ of $\cA$ we have 
\begin{enumerate}
  \item[$1)$] $\im(\I-\phi)$ does not contain any nonzero idempotents of $\cA$. 
  \item[$2)$] $\I-\phi$ maps every $K$-subspace of $\cA$ to a MS of $\cA$.
\end{enumerate}   

In particular, both the LFED conjecture \ref{LFED-Conj} and the LNED conjecture \ref{LNED-Conj} hold for the $K$-$\cE$-derivation $\delta=\I-\phi$ of $\cA$. 
\end{corol}
  
Note that by \cite[Corollary 5.5]{Open-LFNED} both the LFED and LNED conjectures also   
hold for all $K$-$\cE$-derivations associated with finite order $K$-algebra
automorphisms of a commutative $K$-algebra. 
But, for the most of other $K$-algebras these two conjectures are 
still open for this special family of $K$-$\cE$-derivations.

\renewcommand{\theequation}{\thesection.\arabic{equation}}
\renewcommand{\therema}{\thesection.\arabic{rema}}
\setcounter{equation}{0}
\setcounter{rema}{0}

\section{\bf Proof of Theorem \ref{FinExt}}

In this section we   
give a proof for Theorem \ref{FinExt}. We divide the proof into two cases, 
one for the derivation case in subsection \ref{S3.1} and 
the other for the automorphism case  
in subsection \ref{S3.2}.

{\it Throughout this section $R$ stands for a unital ring (not necessarily commutative) 
and $\cA$ for a unital $R$-algebra.}

\subsection{The Derivation Case}\label{S3.1}
We fix an $R$-derivation $D$ of $\cA$ and recall first the construction of the following so-called generalized polynomial algebra over $R$ associated with $D$.

Let $\cA[X; D]$ be the set of all the (generalized) polynomials $f(X)$ 
of the form 
$ 
f(x)=\sum_{i=0}^d a_i X^i
$   
with $d\ge 0$ and $a_i\in \cA$ $(0\le i\le d)$.
Then  $\cA[X; D]$ with the obvious addition and the left scalar multiplication 
forms a left $R$-module. We define a multiplication for $\cA[X; D]$ 
by setting first for all $a\in \cA$ 
\begin{align} 
Xa =aX+Da, \label{[Xa]D}, \\
\intertext{ or equivalently } 
\ad_X (a) =Da, \label{[Xa]D-2}
\end{align}
and then extend it to the product of two arbitrary 
elements of $\cA$ by using the associativity and 
the distribution laws.  

With the operations defined above $\cA[X; D]$ becomes 
a unital $R$-algebra, which contains $\cA$ as an $R$-subalgebra 
and $D=\ad_X \big|_\cA$. In other words,   
the derivation $D$ of $\cA$ becomes inner after $\cA$  
is extended to the $R$-algebra $\cA[X; D]$. 

What we need to show next is that under the integral assumption on $D$ 
we may replace $\cA[X; D]$ by one of its 
quotient $R$-algebras that is finitely generated 
as an $\cA$-module. We begin with the following three lemmas.

\begin{lemma}\label{[Xna]}
For all $n\ge 0$ and $a\in \cA$, we have in $\cA[X, D]$.
\begin{align}\label{[Xna]-eq1}
[X^n, a]\!:=X^na-aX^n=\sum_{i=1}^n \binom ni D^i(a)X^{n-i}.
\end{align}
\end{lemma}

%

\pf First, it is well-known and also easy to check inductively that the following equation holds: 
$$
[X^n, a]=\sum_{i=1}^n \binom ni (\ad_X)^i(a) X^{n-i}.
$$
Then by Eq.\,(\ref{[Xa]D-2}) it is also easy to see 
inductively that $(\ad_X)^i(a)=D^i(a)$ for all $i\ge 0$, from which and the equation above Eq.\,(\ref{[Xna]-eq1}) follows.
\epfv
 

Now, for every $f(X)=\sum_{i=0}^d a_i X^i$ in $\cA[X; D]$, 
we call $a_0$ the {\it constant term} of $f(X)$ 
and denote it by $[0]f(X)$. Since the expression 
$\sum_{i=0}^d a_i X^i$ for $f(X)$ is unique (by definition of 
$\cA[X; D]$), we see that $[0]f(X)$ is well defined.

\begin{lemma} \label{Lma-3.2}
Let $q(X)\in R[X]$ and $(q(X))$ the principal (two-sided) ideal of 
$\cA[X; D]$ generated by $q(X)$. Then for all $f(X)\in (q(X))$, 
we have 
$ 
[0]f(X)\in  \cA \im  q(D),  
$ 
where $\im q(D)\!:=q(D)(\cA)$, i.e., the image of 
the map $q(D):\cA\to \cA$.    
\end{lemma}

\pf By the linearity we may assume $f(X)=a h(X)$ with 
$h(X)=X^m q(X) b X^k $ for some $m, k\ge 0$ 
and $a, b\in \cA$. Therefore it suffices to show $[0]h(X)\in \im q(D)$.

If $k\ge 1$, then it is easy to see by Eq.\,(\ref{[Xa]D}) that 
$[0]h(X)=0$, whence $[0]h(X) \in \im q(D)$.
So assume $k=0$, i.e., $h(X)=X^m q(X)b$ and write 
$q(X)=\sum_{i=0}^d \alpha_i X^i$ for some $d\ge 0$ and 
$\alpha_i\in R$ $(0\le i\le d)$.   

Since $D$ is an $R$-derivation of $\cA$, all $\alpha\in R$ commute with $X$,     
for $\ad_X(\alpha)=D \alpha=0$. Then by Lemma \ref{[Xna]} we have   
\begin{align*}
 q(X)b&=bq(X)+[q(X), b]=bq(X)+[\alpha_0, b]+ \sum_{i=1}^d \alpha_i [X^i, b]\\
&=b q(X)+ [\alpha_0, b] +\sum_{i=1}^d \sum_{j=1}^i \binom ij \alpha_i D^j(b)X^{i-j},   
\end{align*}
from which we get   
\begin{align}
[0] \big( q(X)b \big) &=b \alpha_0 +[\alpha_0, b] + \sum_{i=1}^d \alpha_i D^i(b)
\label{[Xna]-peq2} \\ 
&=\alpha_0b+\sum_{i=1}^d \alpha_i D^i(b)\nno \\&
=\sum_{i=0}^d \alpha_i D^i(b) 
=q(D)(b).\nno
\end{align}

If $m=0$, then $[0]h(X)=[0] \big(q(X) b \big )=q(D)(b)\in \im q(D)$. 
So we assume $m\ge1$.

Since $[X, \alpha]=0$ for all $\alpha\in R$,  
$X$ and $q(X)$ commute. Then 
\begin{align*}
[0]h(X)&=[0] \big( X^m q(X) b \big )=[0] \big( q(X) (X^m b) \big ) \\
&=[0] \big(q(X)b X^m + q(X) [X^m, b]\big)  
\intertext{Applying Lemma \ref{[Xna]}:} 
&=[0] \big(q(X)b X^m + q(X) \sum_{i=1}^m D^i(b)X^{m-i} \big ) \\
&=[0] \big( q(X) D^m(b))  
\intertext{Applying Eq.\,(\ref{[Xna]-peq2}) with $b$ replaced by $D^m(b)$:} 
&=q(D)(D^m(b)).  
\end{align*}
Hence $[0]h(X)\in \im q(D)$, as desired. 
\epfv


\begin{lemma}\label{InjLma-D}
Let $p(t)\in R[t]$ such that $p(D)=0$. 
Then the composition $\cA\to \cA[X; D] \to \cA[X; D]/\big(p(X)\big)$ 
is injective.  
\end{lemma}

\pf Since $p(D)=0$, we have $\im p(D)=0$, whence 
$\cA\im p(D)=0$. Then by Lemma \ref{Lma-3.2} with 
$q(X)=p(X)$ and $f(X)=a$ $(a\in \cA)$ we have 
$\cA \cap (p(X))=0$, from which the lemma follows.
\epfv

Now we are ready to prove Theorem \ref{FinExt} on integral derivations. \\
%

\underline{\bf Proof of Theorem \ref{FinExt}, Part I:} \, 
Let $D$ and $p(t)$ be as in the theorem, and   
set $\cB\!:=\cA[X; D]/\big( p(X)\big)$. 
Then $\cB$ is an $R$-algebra and by Lemma \ref{InjLma-D},  
$\cB$ contains $\cA$ as an $R$-subalgebra. 

To show statement $1)$, denote by $\bar X$ the image of $X$ 
in the quotient algebra $\cB$. Since $p(t)$ is monic, 
it is easy to see that $\cB$ as a left $\cA$-module is 
(finitely) generated by $\bar X^i$ $(0\le i\le d-1)$, 
where $d=\deg p(t)$. Then by Eq.\,(\ref{[Xna]-eq1}) it is easy to verify that 
$\cB$ as a right $\cA$-module is also 
(finitely) generated by $\bar X^i$ 
$(0\le i\le d-1)$.


To show statement $2)$, 
note first that the inner derivation $ad_X$ of $\cA[X, D]$ obviously  
preserves the principal ideal 
$\big(p(X)\big)$. 
Therefore, $ad_X$ induces an inner derivation on the quotient algebra 
$\cB$. Letting $u=\bar X \in \cB$ and by Eq.\,(\ref{[Xa]D-2}) 
we see statement $2)$ also follows. 
\epfv

We end this subsection with the following 
interesting bi-product of Lemma \ref{InjLma-D}.

\begin{corol}\label{PrimeCorol} 
Assume further that $D\ne 0$ and $\cA$ is a simple algebra, i.e., 
the only two-sided ideals of $\cA$ are $0$ and $\cA$ itself. 
Then  $\cA\im D=\cA$.  
\end{corol} 

\pf Assume otherwise, i.e., $\cA\im D\ne \cA$. 
Then by Lemma \ref{Lma-3.2} with $f(X)=q(X)=X$ we have   
$$
(X)\cap \cA \subseteq \cA\im D\ne \cA.
$$ 
Therefore $(X)\cap \cA$ is a proper ideal of $\cA$.  
Since $\cA$ is simple, we have $(X)\cap \cA=\{0\}$.
On the other hand, by Eq.\,(\ref{[Xa]D}) we have 
$\im D\subseteq (X)\cap \cA$, whence $\im D=0$, 
i.e., $D=0$. Contradiction.  
\epfv 

One remark on the corollary above is as follows. 

\begin{rmk}
For all $f(X)\in \cA[X,D]$ write $f(X)=\sum_{i=0}^n X^i c_i$ with $c_i\in \cA$, 
and set $\{0\}f(X)\!:=c_0$. Then by studying $\{0\}f(X)$ instead of $[0]f(X)$ 
it can be shown that the equation $(\im D)\cA=\cA$ also holds in 
Corollary \ref{PrimeCorol}. 
\end{rmk}

\vspace{3mm}

\subsection{The Automorphism Case}\label{S3.2}  

Throughout this subsection we let $R$ and $\cA$ 
be as in Theorem \ref{FinExt}, and $\phi$ 
an  $R$-algebra automorphism of $\cA$, 
and $p(t)$ a monic polynomial in $R[t]$ 
such that $p(\phi)=0$ and $p(0)$ is a unit of $R$.  
 
Let $X$ be a free variable. For any set $S$ we denote by 
$S[X^{-1}, X]$ the set of elements of the form 
$f(X)\!:=\sum_{i=m}^n \alpha_i X^i$ for some $m, n\in \bZ$ and 
$\alpha_i \in S$ $(m\le i\le n)$. Formally, $S[X^{-1}, X]$ is 
just the set of all ``Laurent polynomials" in $X$ with 
the coefficients (appearing on the left) in $S$.

Next we define an associative $R$-algebra $\cA[X^{-1}, X;\phi]$ 
as follows.

First, $\cA[X^{-1}, X; \phi]$ as a set is equal to the set 
$\cA[X^{-1}, X]$.  With the obvious left scalar multiplication 
by elements of $R$ and addition $\cA[X^{-1}, X; \phi]$ 
forms a left $R$-module.  
To make it an $R$-algebra we define the multiplication for 
$\cA[X^{-1}, X; \phi]$ to be the unique associative multiplication 
of $\cA[X^{-1}, X; \phi]$ such that  
$XX^{-1}=X^{-1}X=1_\cA$ and for every $a\in \cA$, 
\begin{align} 
Xa=\phi(a)X,   \label{[Xa]P} \\
\intertext{ or equivalently, } 
 XaX^{-1}=\phi(a). \label{[Xa]P-2}
\end{align}

Then it is easy to see that $\cA[X^{-1}, X; \phi]$ contains $\cA$ 
as an $R$-subalgebra, and by Eq.\,(\ref{[Xa]P-2}) 
we also have $\phi=\Ad_X\,|_{{}_\cA}$.

What we need to show next is that 
we may replace $\cA[X^{-1}, X; \phi]$ by one of its 
quotient algebras $\cB$ which contains $\cA$ and 
is finitely generated as a left $\cA$-module, 
and also as a right $\cA$-module. 

To do so, we need to fix the following notation. 
For every  $f(X)\in \cA[X^{-1}, X; \phi]$, we first write it uniquely as 
$f(X)=\sum_{i=m}^n \alpha_i X^i$ with $\alpha_i\in\cA$ $(m\le i\le n)$, 
and then set $f[1]\!:=\sum_{i=m}^n \alpha_i$. 

One remark on the ``value" $f[1]$ defined above is that,  
if $f(X)=\sum_{i=k}^\ell X^i\beta_i$ with $\beta_i\in\cA$ $(k\le i\le \ell)$, 
then $f[1]$ may not be equal to 
$\sum_{i=0}^n \beta_i$. In other words, 
in order to ``evaluate $f(X)$ 
at $X=1$", we need first to write it as a Laurent polynomial in $X$ 
with all the coefficients appearing on the left of $X^i$'s.

\begin{lemma}\label{Lma-3.2.1}
Let $f(X)\in R [X^{-1}, X]$ and 
$\big(f(X)\big)$ the principal (two-sided) ideal of 
$\cA[X^{-1}, X; \phi]$ generated by $f(X)$. 
Then for every $h(X)\in \big(f(X)\big)$ we have  
\begin{align}
h[1]\in \cA \im f(\phi),  \label{Lma-3.2.1-eq1}
\end{align}
where $\im f(\phi)\!:=f(\phi)(\cA)$. 
\end{lemma}

\pf Write $f(X)=\sum_{i=m}^n \alpha_iX^i$ for some $m, n\in \bZ$ and 
$\alpha_i$'s in $R$. By the linearity we may assume 
$h(X)=bX^j f(X)cX^k$ for some $j, k\in \bZ$ and $b, c\in \cA$. Since $\phi$ 
is an $R$-algebra automorphism of $\cA$, we have $\phi(r1_\cA)=r1_\cA$. Then by 
Eq.\,(\ref{[Xa]P}) we see that all elements of $R$ commute with $X$ in 
$\cA[X^{-1}, X; \phi]$. By this fact and Eq.\,(\ref{[Xa]P-2}) we consider   
\begin{align*}
h(X)&=bX^j f(X)cX^k = b \sum_{i=m}^n  (X^j \alpha_i) (X^ic) X^k \\ 
&= b \sum_{i=m}^n  (X^j \alpha_i) (X^icX^{-i}) X^{i+k} 
 =b  \sum_{i=m}^n  \alpha_i X^j  \phi^i(c) X^{i+k}\\&
=b\sum_{i=m}^n \alpha_i \phi^j (  \phi^i(c)) X^{j+i+k} 
 =b \big(\sum_{i=m}^n  \alpha_i \phi^i(\phi^j(c))\big) X^{j+i+k}.
\end{align*}
Therefore we have 
$$
h[1]= b \big(\sum_{i=m}^n  \alpha_i \phi^i(\phi^j(c))\big)=b f(\phi)(\phi^j(c)).
$$
Hence the lemma follows. 
\epfv

\begin{lemma}  \label{Lma-3.2.2}
For each $f(X)\in R [X^{-1}, X]$ such that $f(\phi)=0$, 
the composition $\cA \to \cA[X^{-1}, X; \phi] 
\to \cA[X^{-1}, X; \phi]/\big(f(X)\big)$ is injective.  
\end{lemma}

\pf Since $f(\phi)=0$, $\cA \im f(\phi)=\{0\}$. 
Then for each $a\in \cA$, by Lemma \ref{Lma-3.2.1} we have that $a\in \big(f(X)\big)$, 
if and only if $a=0$. 
Therefore $\cA \cap \big(f(X)\big)=\{0\}$, whence the lemma follows. 
\epfv

Now we are ready to show the automorphism case of Theorem \ref{FinExt}.\\ 


\underline{\bf Proof of Theorem \ref{FinExt}, Part II:}\, 
 Let $\phi$ and $p(t)\in R[t]$ be as in the theorem, and set   
 $\cB\!:=\cA[X^{-1}, X; \phi]/ \big( p(X) \big)$. 
Then by Lemma \ref{Lma-3.2.2} $\cB$ contains $\cA$ 
as an $R$-subalgebra. 
Since the inner automorphism 
$\Ad_X$ of $\cA[X^{-1}, X; \phi]$ obviously preserves the principal ideal 
$\big( p(X) \big)$, $\Ad_X$ induces an inner automorphism of $\cB$. Let $u$ be the image of $X$ 
in $\cB$. Then  $p(u)=0$ and by Eq.\,(\ref{[Xa]P-2}) $\phi=\Ad_u\big|_\cA$.   

Therefore, it remains only to show that $\cB$ as a left $\cA$-module as well as 
a right $\cA$-module is finitely generated. 
Write $p(t)=t^d+\sum_{i=0}^{d-1} \alpha_i t^i$ 
for some $d\ge 0$ and $\alpha_i\in R$ $(0\le i\le d)$.
Then $\bar X^d=-\sum_{i=0}^{d-1} \alpha_i \bar X^i$ in $\cB$, 
whence for each $m\ge d$, $\bar X^m$ lies in the left 
$\cA$-submodule generated by $\bar X^i$ $(0\le i\le d-1)$. 

On the other hand, since by assumption $\alpha_0=p(0)$ is a unit in $R$, and  
$\alpha_0^{-1}\bar X^{-d}p(\bar X)=0$ in $\cB$, we have 
$$
\bar X^{-d}=-\alpha_0^{-1}-\sum_{i=1}^{d-1} \alpha_0^{-1}\alpha_i \bar X^{-d+i}.
$$
Therefore, for each $n\le -d$, 
$\bar X^n$ lies in the $\cA$-submodule generated by $\bar X^{-i}$ 
$(0\le i\le d-1)$. Hence $\cB$ as a left $\cA$-module 
is (finitely) generated by $u^i=\bar X^{i}$ $(-d+1\le i\le d-1)$.

Furthermore, by Eq.\,(\ref{[Xa]P}) it is easy to see inductively  
that $\bar X^k r=\phi^k(r) \bar X^k$ for all $a\in \cA$ and $k\in \bZ$, 
from which we see that $\cA$ as a right $\cA$-module is also (finitely) 
generated by $\bar X^{i}$ $(-d+1\le i\le d-1)$.    
\epfv

\end{document}